\documentclass{article}
\usepackage{authblk}

\newtheorem{fed}{\textbf{Definition}}[section]
\newtheorem{thm}[fed]{\textbf{Theorem}}
\newtheorem{lemma}[fed]{\textbf{Lemma}}

\newtheorem{cor}[fed]{\textbf{Corollary}}
\usepackage{amssymb,bbm,graphicx,epsfig,psfrag,epic,eepic,latexsym}
\usepackage{amsmath}
\usepackage{mathrsfs}
\usepackage[dvips]{color}

\begin{document}
\title{A Hamiltonian version of a result of Gromoll and Grove}
\author[1]{Urs Frauenfelder}
\author[2]{Christian Lange}
\author[3]{Stefan Suhr}
\affil[1]{Institut f\"ur Mathematik, Universit\"at Augsburg, Universit\"atsstrasse 14, 86159 Augsburg, Germany}
\affil[2]{Mathematisches Institut der Universit\"at zu K\"oln, Weyertal 86-90, 50931 Cologne, Germany}
\affil[3]{DMA, ENS Paris and Universit\'e Paris Dauphine, 45 rue d'Ulm, 75005 Paris, France}
\maketitle
\abstract{The theorem that if all geodesics of a Riemannian two-sphere are closed they are also simple closed is generalized to real Hamiltonian structures 
on $\mathbb{R}P^3$. For reversible Finsler $2$-spheres all of whose geodesics are closed this implies that the lengths of all geodesics coincide.}
\section{Introduction}
In \cite{poincare} Poincar\'e advertised the study of geodesics as a first step to penetrate in
the exploration of Hamiltonian systems like the ones appearing in celestial mechanics. A fascinating
object in the study of geodesics are manifolds all whose geodesics are closed, see \cite{besse}. This leads
to the question how much results on closed geodesics are Riemannian phenomena or belong to the realm of
Hamiltonian dynamics. For example in \cite{frauenfelder-labrousse-schlenk} it was shown how the classical
result of Bott-Samelson fits into a contact geometric set-up. 

The theorem of Gromoll and Grove \cite{gromoll-grove} asserts that if a Riemannian metric on the two dimensional sphere $S^2$ has the property that all geodesics are closed, then each geodesic is simple closed
and therefore all geodesics have a \emph{common} minimal period. 

The geodesic equation is a second order ODE. Each second order ODE can be transformed into a first order ODE. This allows one to interpret the geodesic equation as a first order Hamiltonian flow equation on the unit cotangent bundle of the underlying manifold. In the case of $S^2$ the unit cotangent bundle $S^*S^2$ is
diffeomorphic to real projective space $\mathbb{R}P^3$. In physics the underlying manifold is referred to
as the configuration space while its cotangent space is denoted as the phase space. 

The phase space is a symplectic manifold and the symplectic geometer is interested to apply arbitrary symplectic transformations to this symplectic manifold, since under these transformations the symplectic dynamics transforms canonically and can under some lucky instances brought to a form which makes 
it more tractible to explicit solutions. But under arbitrary symplectic transformations 
the bundle structure of the cotangent bundle is not preserved. Therefore the question if a closed characteristic of the Hamiltonian flow
is simple or not does not make sense to the symplectic geometer. However, in the case that all characteristics
are closed, if the Hamiltonian flow admits a \emph{common} minimal period, is a question even a symplectic geometer can understand. 

Unfortunately, in general this is bloody wrong. Here is an example. Let $p$ and $q$ be two relatively prime
integers. Identify the three dimensional sphere $S^3$ with $\{(z_1,z_2) \in \mathbb{C}^2: |z_1|^2+|z_2|^2=1\}$ and consider the circle
action on $S^3$ given by
$$e^{it}_*(z_1,z_2)=(e^{ipt}z_1,e^{iqt}z_2).$$
This circle action commutes with the antipodal involution $z \mapsto -z$ on $S^3$ and hence induces a circle action on $\mathbb{R}P^3=S^3/\mathbb{Z}_2$. The minimal period of the orbit through $[(1,0)]$ is $\tfrac{\pi}{p}$, the minimal
period of the orbit through $[(0,1)]$ is $\tfrac{\pi}{q}$, while the minimal period of an orbit through a point
$[(z_1,z_2)] \in \mathbb{R}P^3$ with $z_1 \neq 0$ and $z_2\neq 0$ is $\pi$ in the case $p$ and $q$ are both odd and
$2 \pi$ if one of them is even. It is worth noting that this example corresponds to the Katok examples as was pointed out by Harris and Paternain \cite{harris-paternain}. This means that the Gromoll-Grove theorem
already fails for nonreversible Finsler metrics.

Geodesic flows of a Riemannian or a reversible Finsler metric have the following properties which distinguish them from other Hamiltonian flows
\begin{description}
 \item[(i)] The geodesic flow is invariant under time reversal, i.e., a geodesic traversed backwards is again
a geodesic.
 \item[(ii)] The geodesic flow is a Reeb flow, i.e., it is the flow of a Reeb vector field of a contact structure.
\end{description}
In this note we explain how under a suitable generalization of property (i) for geodesics the Gromoll-Grove result
about the common minimal period generalizes to hold for general Hamiltonian systems all of whose trajectories are closed on energy hypersurfaces
having the same topology as the ones studied by Gromoll and Grove, namely $\mathbb{R}P^3$. Here we do not need to assume that our Hamiltonian flow is a Reeb flow. Instead of that we explain invoking a deep result
of Epstein \cite{epstein} that in this set-up the Hamiltonian flow will be stable. 

The proof of Gromoll and Grove hinges on the Theorem of Lusternik and Schnirelmann on the existence of
three simple closed geodesics and the theory of Seifert on three dimensional fibred spaces. As explained above there is no
analogon of the Theorem of Lusternik and Schnirelmann for Hamiltonian flows because the notion of simpleness
does not generalize to the symplectic set-up. However, the tangent curve of a smooth regular simple curve (i.e. an embedding) $S^1\to S^2$ 
is a noncontractible curve in $TS^2\setminus S^2$. We will show that the noncontractibility of all orbits follows without further assumptions. 
The analogue in the Riemannian case is proven in \cite{lange}. Here we emply the same method.
%under the assumption of at least three noncontractible orbits of a real 
%Hamiltonian structure (Definition \ref{realham}) on $\mathbb{R}P^3$ all orbits have a common minimal period. Further we conjecture that without this assumption 
%the theorem fails in general, i.e. one cannot do without an analogue of the Theorem of Lustenik and Schnirelman. If this indeed happens one can 
%wonder whether the Theorem of Lusternik and Schnirelman holds for reversible Finsler metrics. In other words: Do the conjectural counterexamples 
%admit a property which distinguishes them from reversible Finsler metrics. 
\\ \\
\emph{Acknowledgements: } Urs Frauenfelder would like to thank the Ecole Normale Sup\'erieure for its
hospitality during the visit while this research was carried out.  

 The research leading to these results has received funding from the
European Research Council under the European Union's Seventh Framework
Programme (FP/2007-2013) / ERC Grant Agreement 307062.

\section{Real Hamiltonian manifolds}
Assume that $\Sigma=\Sigma^{2n+1}$ is a closed oriented manifold of dimension $2n+1$. A \emph{Hamiltonian structure} on $\Sigma$ is a two form
$$\omega \in \Omega^2(\Sigma)$$
satisfying the following two conditions
\begin{description}
 \item[(i)] $\omega$ is closed, i.e., $d\omega=0$,
 \item[(ii)] $\mathrm{ker}(\omega)$ defines a one dimensional distribution on $T\Sigma$.
\end{description}
We refer to the pair $(\Sigma,\omega)$ is a \emph{Hamiltonian manifold}. A Hamiltonian manifold
is the odd dimensional analogon of a symplectic manifold. In the study of Hamiltonian dynamics Hamiltonian manifolds naturally arise as energy hypersurfaces, i.e., level sets of a smooth function on a symplectic manifold. The characteristic foliation, namely the leaves of the kernel of $\omega$, correspond to the trajectories of the Hamiltonian flow on a fixed energy level. People studying Hamiltonian manifolds are relaxed people, they do not care about the parametrization of the leaves but only about the foliation itself. 

Because we assume our manifold $\Sigma$ to be oriented we actually obtain an orientation of the leaves. Indeed, note that for each $x \in \Sigma$ the two form $\omega_x$ is nondegenerate on $T_x \Sigma/\mathrm{ker}(\omega_x)$, hence is a symplectic form on this vector space and therefore defines an orientation on it. 
\begin{fed}\label{realham}
A \emph{real Hamiltonian manifold} is a triple $(\Sigma, \omega, \rho)$, where
$(\Sigma,\omega)$ is a Hamiltonian manifold and $\rho \in \mathrm{Diff}(\Sigma)$ is a smooth involution on $\Sigma$, i.e., $\rho^2=\mathrm{id}_\Sigma$ satisfying the following two conditions
\begin{description}
 \item[(i)] $\rho^* \omega=-\omega$,
 \item[(ii)] $\rho$ reverses the orientation on $\mathrm{ker}(\omega)$.
\end{description}
\end{fed}
If $\mathscr{L}$ denotes the set of leaves of the foliation on $\Sigma$ induced from $\mathrm{ker}(\omega)$,
then the first condition guarantees that $\rho$ induces an involution on the set of leaves
$$\rho_* \colon \mathscr{L} \to \mathscr{L}.$$
Here is an example of a real Hamiltonian manifold. For a closed Riemannian manifold $(N,g)$ let 
$$\Sigma=S^* N=\{(q,p) \in T^*N: ||p||_{g_q}=1\}$$
be its unit cotangent bundle with Hamiltonian structure obtained by restricting the differential of the Liouville one-from to $S^* N$. Then a real structure is given by
$$\rho \colon S^* N \to S^*N, \quad (q,p) \mapsto (q,-p).$$
Note that the same construction works for reversible Finsler metrics but fails for nonreversible Finsler metrics since in this case the unit cotangent bundle is not invariant under $\rho$ anymore. In this example the leaves of the foliation correspond to oriented geodesics and the involution reverses the orientation of the geodesic. In particular, no leaf is fixed. This is a more general phenomenon as the following Lemma shows. 
\begin{lemma}\label{free}
Suppose that $(\Sigma, \omega,\rho)$ is a real Hamiltonian manifold with the property that 
$$\mathrm{Fix}(\rho)=\emptyset,$$
i.e., there are no fixed points of the involution $\rho$ on $\Sigma$. Then the induced involution
$\rho_*$ on the set of leaves $\mathscr{L}$ is free. 
\end{lemma}
\textbf{Proof: } We argue by contradiction and assume that there exists a leaf $L \in \mathscr{L}$ such that
$$\rho_*(L)=L.$$
Hence if $x \in L$ it holds that
$$\rho(x) \in L.$$
Because $\rho$ has no fixed points by assumption we must have
$$x \neq \rho(x).$$
The leaf $L$ is either diffeomorphic to a circle or to the real line. Hence we can choose
$$I \subset L$$
diffeomorphic to a closed interval such that
$$\partial I=\{x, \rho(x)\}.$$
Note that since $\rho$ is an involution it interchanges the boundary points of $I$. Because $\rho$ reverses the orientation of $L$ we conclude that
$$\rho(I)=I.$$
However, because $\rho$ interchanges the boundary points of the closed interval $I$ we conclude by the intermediate value theorem that there exists $y \in I$ such that
$$\rho(y)=y$$
contradicting the assumption that $\mathrm{Fix}(\rho)=\emptyset$. This contradiction concludes the proof of the lemma. \hfill $\square$
\section{Stability}

Assume that $(\Sigma,\omega)$ is a Hamiltonian manifold of dimension $2n+1$. A one-form 
$\lambda \in \Omega^1(\Sigma)$ is called a \emph{stabilizing one-form}, see \cite{cieliebak-frauenfelder-paternain, cieliebak-volkov}, if it satisfies the following two conditions
\begin{description}
 \item[(i)] $\lambda \wedge \omega^n>0$, i.e., $\lambda \wedge \omega^n$ is a volume form on
 the oriented manifold $\Sigma$,
 \item[(ii)] $\ker( \omega) \subset \ker(d\lambda)$.
\end{description}
We abbreviate by
$$\Lambda(\Sigma,\omega) \subset \Omega^1(\Sigma)$$
the maybe empty subset of all stabilizing one-forms of $(\Sigma,\omega)$. Note that
$\Lambda(\Sigma,\omega)$ is a cone, i.e., if $\lambda_1, \lambda_2 \in \Lambda(\Sigma,\omega)$
and $r_1, r_2>0$, then $r_1 \lambda_1+r_2 \lambda_2 \in \Lambda(\Sigma,\omega)$. We refer to
$\Lambda(\Sigma,\omega)$ as the \emph{stable cone} of $(\Sigma,\omega)$. 
\begin{fed}
A Hamiltonian manifold $(\Sigma,\omega)$ is called \emph{stable}, if $\Lambda(\Sigma,\omega) \neq \emptyset$, i.e., there exists a stabilizing one-form.
\end{fed}
Now suppose that $(\Sigma,\omega,\rho)$ is a real Hamiltonian manifold. We say that 
a stabilizing one-form $\lambda \in \Lambda(\Sigma,\omega)$ is \emph{real} if it satisfies
$$\rho_* \lambda=-\lambda.$$
To produce real stabilizing one-forms we first show that the cone $\Lambda(\Sigma,\omega)$
is invariant under the involution $\lambda \mapsto -\rho^* \lambda$ of $\Omega^1(\Sigma)$.
\begin{lemma}\label{realization}
Suppose that $\lambda \in \Lambda(\Sigma,\omega)$. Then
$$-\rho^* \lambda \in \Lambda(\Sigma, \omega).$$
\end{lemma}
\textbf{Proof: }From the definition of a real Hamiltonian manifold we see immediately that 
$$(-1)^{n+1} \rho^*(\lambda\wedge \omega^n) >0.$$
Thus we have 
$$(-\rho^*\lambda)\wedge \omega^n=(-1)^{n+1}\rho^*(\lambda\wedge\omega)>0,$$
and the first condition for a stablizing one-form is satisfied by $-\rho^*\lambda$.

Combining $\ker(\rho^*\omega)=\ker \omega$ with $\ker\omega\subset \ker d\lambda$ implies 
$$\ker\omega=\ker \rho^*\omega\subset \ker \rho^*d\lambda=\ker (-\rho^*\lambda).$$
Therefore $-\rho^*\lambda$ is a stabilizing one-form. \hfill $\square$
\\ \\
In view of the cone property of $\Lambda(\Sigma, \omega)$ the lemma enables us to construct
real stabilizing one-forms out of a stabilizing one-form.
\begin{cor}
Suppose that the real Hamiltonian manifold $(\Sigma,\omega,\rho)$ is stable. Then it admits a real
stabilizing one-form. 
\end{cor}
\textbf{Proof: } Let $\lambda$ be a stabilizing one-form. By Lemma~\ref{realization} it follows that
$-\rho^* \lambda$ is again a stabilizing one-form. Because $\Lambda(\Sigma,\omega)$ is a cone it follows that 
$$\lambda_\rho:=\lambda-\rho^* \lambda$$
is again a stabilizing one-form. Because $\rho$ is an involution we have
$$\rho^* \lambda_\rho=-\rho^* \lambda_\rho.$$
Therefore $\lambda_\rho$ is a real stabilizing one-form. \hfill $\square$
\section{Periodic Hamiltonian manifolds}
\begin{fed}
A Hamiltonian manifold $(\Sigma,\omega)$ is called \emph{periodic} if all leaves are circles.
\end{fed}
\begin{thm}\label{stab}
Assume that $(\Sigma,\omega)$ is a $3$-dimensional periodic Hamiltonian manifold. Then it is stable. 
\end{thm}
\textbf{Proof: }For every Hamiltonian structure on an orientable manifold one can orient the characteristic distribution. Therefore by choosing a 
nowhere vanishing vectorfield tangent to the charateristic distribution one finds that the charateristic foliation consists of orbits of a $\mathbb{R}$ 
action on $\Sigma$. 

According to the main theorem in \cite{epstein} if $(\Sigma,\omega)$ is periodic on the $3$-manifold $\Sigma$ then the $\mathbb{R}$ action 
is orbit equivalent to an $S^1$ action. Thus the length of the leaves with respect to any Riemannian metric on $\Sigma$ is bounded. 
Then by the Theorem of Wadsley in the version of \cite{cieliebak-volkov} there exists a stabilizing one-form.\hfill $\square$

\begin{cor}\label{corstab}
Under the assumptions of Theorem~\ref{stab} there exists a smooth circle action on $\Sigma$
without fixed points such that the orbits of the circle action correspond to the leaves of the
characteristic foliation of $(\Sigma,\omega)$.
\end{cor}
\section{The main result and proof}
\begin{thm}\label{main}
Suppose that $(\omega,\rho)$ is a periodic real Hamiltonian structure on $\mathbb{R}P^3$. %with the property that $\mathrm{Fix}(\rho)=\emptyset$ and at least three leaves are noncontractible. 
Then there exists a free circle action on $\mathbb{R}P^3$ 
with the property that the leaves of $(\mathbb{R}P^3, \omega)$ correspond to the orbits of the circle action. 
In particular, all leaves are noncontractible.
\end{thm}

\begin{cor}
Let $(S^2,F)$ be a reversible Finsler $2$-sphere such that all geodesics are closed. Then all geodesics give rise to noncontractible orbits on the unit tangent bundle of $S^2$ and have the same length.
\end{cor}

%\begin{rem}
%By Lemma~\ref{free} the induced involution $\rho_*$ on the space of leaves is free. Therefore if there
%are at least three noncontractible leaves there have to be actually at least four noncontractible leaves. 
%\end{rem}
%
%\begin{rem}
%The theorem of Lusternik-Schnirelmann \cite{lusternik-schnirelmann} guarantees for every Riemannian
%metric on $S^2$ at least three simple closed geodesics. Each simple closed geodesic gives rise to two
%noncontractible periodic orbits on the unit cotangent bundle $S^*S^2=\mathbb{R}P^3$. The fact that there 
%are two of them is because each geodesic has two different orientations. Therefore the theorem of 
%Lusternik-Schnirelmann provides the existence of at least six noncontractible closed leaves. 
%\end{rem}
%\section{Proof of the main result}

In his fundamental monograph \cite{seifert} Seifert studies effective circle actions on closed three dimensional manifolds
with no fixed points. For a modern account of Seifert's theory see for example the paper by Scott
\cite{scott}. If the three dimensional manifold is orientable the quotient orbifold is an orientable surface with
finitely many cone points. If $S$ is a surface and $p_1, \cdots p_\ell$ is a finite collection of positive integers
bigger than two we abbreviate by $S(p_1, \cdots, p_\ell)$ the surface with $\ell$ cone points of order
$p_1, \cdots, p_\ell$. 

\begin{thm}[Seifert]\label{seifert}
Suppose that $S^1$ acts smoothly on $\mathbb{R}P^3$ without fixed points. Then the quotient orbifold
$\mathbb{R}P^3/S^1$ belongs to one of the following examples
$$S^2, \quad S^2(p), \quad S^2(p,q).$$
In the third case it holds that
$$(p,q) \in \{1,2\}$$
i.e., $p$ and $q$ are either relatively prime or decompose as $p=2p'$ and $q=2q'$ with
$p'$ and $q'$ relatively prime. 
\end{thm}
\textbf{Proof: } For the fibration $S^1 \to \mathbb{R}P^3 \to \mathbb{R}P^3/S^1$ we obtain an exact sequence of homotopy groups
$$\pi_1(S^1) \rightarrow \pi_1(\mathbb{R}P^3) \rightarrow \pi_1(\mathbb{R}P^3/S^1)
\to \pi_0(S^1)$$
see \cite[Lemma 3.2]{scott}. Hence we have an exact sequence
$$\mathbb{Z} \rightarrow \mathbb{Z}_2 \rightarrow \pi_1(\mathbb{R}P^3/S^1) \to \{1\}$$
implying that the fundamental group of the orbifold $\mathbb{R}P^3/S^1$ is either trivial or $\mathbb{Z}_2$.
Closed two dimensional orbifolds of trivial fundamental group or fundamental group $\mathbb{Z}_2$
are the three examples in the list above as we explain in the appendix. This finishes the proof of the Theorem. \hfill $\square$
\begin{cor}\label{corseif}
Suppose that $(\omega,\rho)$ is a real Hamiltonian structure on $\mathbb{R}P^3$ with the property that
$\mathrm{Fix}(\rho)=\emptyset$. Then the space of leaves $\mathscr{L}$  is diffeomorphic to $S^2$ and all orbits are noncontractible. 
%$$S^2, \quad S^2(2,2).$$
\end{cor}
Theorem \ref{main} follows directly from the corollary. \\ \\
\textbf{Proof: } By Corollary~\ref{corstab} there exists a circle action on $\mathbb{R}P^3$ without
fixed points such that the space of leaves is
$$\mathscr{L}=\mathbb{R}P^3/S^1.$$
Because the involution $\rho$ has no fixed points by assumption Lemma~\ref{free} tells us that 
the induced involution $\rho_*$ on the space of leaves $\mathscr{L}$ is free. In particular, the number of cone
points of a given multiplicity is even. Therefore by Theorem~\ref{seifert} the space of leaves belongs to one of the following orbifolds
$$S^2, \quad S^2(p,p)$$
where in the later case we have
$$(p,p) \in \{1,2\}$$
implying that 
$$p=2.$$
 Lemma 3.4 in \cite{lange} shows that every effective $S^1$-action without fixed points on a lens space $L(r,1)$ with orbit space $S^2(2,2)$ implies $r\ge 4$. Here the lens space $L(r,1)$ is defined 
as $S^3/\mathbb{Z}_r$ where the $\mathbb{Z}_r$-action is generated by multiplication with $e^{\frac{2\pi}{r}i}$. This excludes the orbifold $S^2(2,2)$ as the orbit space since 
$\mathbb{R}P^3 \cong L(2,1)$. Using the homotopy exact sequence
$$\pi_1(S^1)\to \pi_1(\mathbb{R}P^3)\to \pi_1(S^2)=0$$
we see that all orbits are noncontractible.
This finishes the proof of the Corollary. \hfill $\square$
%\\ \\
%\textbf{Proof of Theorem~\ref{main}: }By Corollary~\ref{corseif} the space of leaves is either
%$S^2$ or $S^2(2,2)$. In the later case all regular leaves are homotopic to the double covering of the
%exceptional ones and hence because the fundamental group of $\mathbb{R}P^3$ is $\mathbb{Z}_2$
%have to be contractible. Therefore there are at most two noncontractible leaves, namely the
%exceptional ones contradicting the assumption of the theorem. Consequently the space of leaves
%is $S^2$. Because $S^2$ is a manifold and has therefore no cone points the action of
%$S^1$ on $\mathbb{R}P^3$ is free. In particular, all leaves are homotopic to each other and because
%three of them are noncontractible all of them are noncontractible. This finishes the proof of the main result.
%\hfill $\square$

\appendix

\section{The orbifold fundamental group}

Assume that $\Sigma$ is an orientable closed 3-dimensional manifold on which the circle $S^1$ acts 
effectively without fixed points. Then the quotient $S:=\Sigma/S^1$ is an orbifold, namely a closed orientable surface with finitely many cone points (say $\ell$ many) of multiplicity $p_i$ for $1 \leq i \leq \ell$. Suppose that the genus of the surface $S$ is $g$. Interpreted as a topological space
$S_{\mathrm{top}}$ and not as orbifold its fundamental group has the presentation
$$\pi_1(S_{\mathrm{top}})=\Big\langle a_1,b_1, \cdots a_g, b_g\big| a_1 b_1 a_1^{-1} b_1^{-1} \cdots a_g b_g a_g^{-1} b_g^{-1} \Big\rangle.$$
Considered as an orbifold $S_{\mathrm{orb}}$ one defines its \emph{orbifold fundamental group} as
\begin{eqnarray*}
\pi_1(S_{\mathrm{orb}})=\qquad \qquad \qquad \qquad \qquad \qquad \qquad \\
\Big\langle a_1,b_1, \cdots a_g, b_g, x_1, \cdots, x_\ell\big| x_1^{p_1}, \cdots, x_\ell^{p_\ell}, a_1 b_1 a_1^{-1} b_1^{-1} \cdots a_g b_g a_g^{-1} b_g^{-1}x_1 \cdots x_p \Big\rangle.
\end{eqnarray*}
The usefulness of the orbifold fundamental group lies in the fact that one has a surjective group homomorphism
$$\pi_1(\Sigma) \twoheadrightarrow \pi_1(S_{\mathrm{orb}})$$
as was proved by Seifert in \cite[Chapter 10]{seifert}. This comes from the homotopy exact sequence
$$\pi_1(S^1) \to \pi_1(\Sigma) \to \pi_1(S_{\mathrm{orb}}) \to \pi_0(S^1),$$
see for example \cite[Chapter 3]{scott}. The existence of this surjective group homomorphism gives us interesting information on the types of exceptional fibers the action of $S^1$ on $\Sigma$ can have. 

Suppose in the following that the fundamental group of $\Sigma$ is finite. In view of the surjective group homomorphism the orbifold fundamental group of $S$ is finite as well. In particular, topologically
$S$ has to be a sphere
$$S_{\mathrm{top}}=S^2$$
and its orbifold fundamental group has the presentation
\begin{eqnarray*}
\pi_1(S_{\mathrm{orb}})=
\Big\langle x_1, \cdots, x_\ell\big| x_1^{p_1}, \cdots, x_\ell^{p_\ell}, x_1 \cdots x_p \Big\rangle.
\end{eqnarray*}
It turns out that this group is finite only if there are less than three cone points or if there are three cone points
the orbifold is of the form
$$S^2(2,2,n), \quad S^2(2,3,3), \quad S^2(2,3,4), \quad S^2(2,3,5)$$
for $n \geq 2$. For a beautiful geometric interpretation why this holds true we invite the reader to
consult directly Chapter\,10 in Seifert's paper \cite{seifert}.
In the case where there is just one cone point, the orbifold is a teardrop, and its fundamental group is trivial.
If there are two cone points the fundamental group is
$$\pi_1(S^2(p,q))=\mathbb{Z}/(p,q)\mathbb{Z}$$
where $(p,q)$ is the greatest common divisor of the integers $p$ and $q$. The fundamental group
of $S^2(2,2,n)$ is the dihedral group of order $2n$ while the fundamental groups $S^2(2,3,3)$, $S^2(2,3,4)$ and
$S^2(2,3,5)$ correspond to the symmetry groups of the platonic solids. In particular, we see that if the orbifold
fundamental group of $S$ is trivial or $\mathbb{Z}_2$, then the orbifold belongs to one of the examples listed
in Theorem~\ref{seifert}.

\end{document}